\documentclass[10pt,preprint]{article}

\usepackage{amsmath, latexsym, amsfonts, amssymb, amsthm}
\usepackage{graphicx, color,hyperref,dsfont,epsfig,caption,wrapfig,subfig}
\usepackage{pstricks}
\usepackage{movie15}
\hypersetup{
    linktoc=page,
    linkcolor=red,          
    citecolor=blue,        
    filecolor=blue,      
    urlcolor=cyan,
    colorlinks=true           
}

\numberwithin{equation}{section}

\newtheorem{theorem}{Theorem}[section]
\newtheorem{lemma}[theorem]{Lemma}
\newtheorem{proposition}[theorem]{Proposition}

\newtheorem{definition}[theorem]{Definition}

\newcounter{conj}
\newtheorem{conjecture}[conj]{Conjecture}

\theoremstyle{definition}

\renewcommand{\tilde}{\widetilde}          
\DeclareMathSymbol{\leqslant}{\mathalpha}{AMSa}{"36} 
\DeclareMathSymbol{\geqslant}{\mathalpha}{AMSa}{"3E} 
\DeclareMathSymbol{\eset}{\mathalpha}{AMSb}{"3F}     
\renewcommand{\leq}{\;\leqslant\;}                   
\renewcommand{\geq}{\;\geqslant\;}                   

\newcommand{\C}{\mathbb{C}}

\newcommand{\R}{\mathbb{R}}
\newcommand{\Z}{\mathbb{Z}}

\newcommand{\E}{\mathds{E}}
\renewcommand{\P}{\mathds{P}}
\newcommand{\ind}{\mathds{1}}

\usepackage{xparse}

\DeclareDocumentCommand \Pmp { m m o} {
\IfNoValueTF{#3}
{P_{#1}^{#2}}
{P_{#1}^{#2}\left(#3\right)}
}

\DeclareDocumentCommand \Emp { m m o} {
\IfNoValueTF{#3}
{E_{#1}^{#2}}
{E_{#1}^{#2}\left[#3\right]}
}

\DeclareDocumentCommand \Pbr { m m m m o } {
\IfNoValueTF{#5}
{P_{#1}^{#2\stackrel{#4}{\rightarrow} #3}}
{P_{#1}^{#2\stackrel{#4}{\rightarrow} #3}\left(#5\right)}
}

\DeclareDocumentCommand \Ebr { m m m m o } {
\IfNoValueTF{#5}
{E_{#1}^{#2\stackrel{#4}{\rightarrow} #3}}
{E_{#1}^{#2\stackrel{#4}{\rightarrow} #3}\left[#5\right]}
}

\def\S{\mathbb{S}}

\def\bi{\begin{itemize}}
\def\ei{\end{itemize}}

\def\bnum{\begin{enumerate}}
\def\enum{\end{enumerate}}
\def\<#1{\langle #1 \rangle}

\title{Lecture notes on Gaussian multiplicative chaos and Liouville Quantum Gravity}

 \author{  R\'emi Rhodes \footnote{Universit{\'e} Paris-Est Marne la Vall\'ee, LAMA, Champs sur Marne, France.} \footnotetext[2]{Partially supported by grant ANR-11-JCJC  CHAMU.},
 Vincent Vargas \footnote{ENS Ulm, DMA, 45 rue d'Ulm,  75005 Paris, France.} }

\date{\vspace{-5ex}}

\begin{document}

\maketitle

\begin{abstract}
The purpose of these notes, based on a course given by the second author at Les Houches summer school, is to explain the probabilistic construction of Polyakov's Liouville quantum gravity using the theory of Gaussian multiplicative chaos. In particular, these notes contain a detailed description of the so-called Liouville measures of the theory and their conjectured relation to the scaling limit of large planar maps properly embedded in the sphere. These notes are rather short and require no prior knowledge on the topic.   
\end{abstract}

\tableofcontents
\noindent
\section{Introduction}

In 1985, Kahane laid the foundations of Gaussian multiplicative chaos theory (GMC, hereafter).  Roughly speaking, GMC is a theory which defines rigorously random measures with the following formal definition
\begin{equation}\label{firstdef}
M_\gamma (dx)=e^{\gamma X(x)} \sigma(dx)
\end{equation}
where $\sigma$ is a Radon measure on some metric space $D$ (equipped with a metric $d$), $\gamma>0$ is a parameter and $X:D\to\R$ is a centered Gaussian field. The definition \eqref{firstdef} should be seen as formal since in the interesting cases the variable $X$ does not live in the space of functions on $D$ but rather in a space of distributions in the sense of Schwartz. In that case, $X(x)$ does not make sense pointwise. Of course, we could make the change of variables $X \to \gamma X$ and absorb the dependence in $\gamma$ in the field $X$ but we will not do so for reasons which will become clear in the sequel. In fact, Kahane's GMC theory is quite general and the metric space $D$ need not be some subspace of $\R^d$; however, motivated by the study of 2d Liouville quantum gravity (LQG, hereafter), we will consider in the sequel the very important subcase where $D$ is some subdomain of $\R^d$, $\sigma$ is  a Radon measure on $D$ and $X$ has a covariance kernel of log-type, namely
\begin{equation}\label{logcov}
K(x,y) := \E[X(x)X(y)]=\ln_+\frac{1}{|x-y|}+g(x,y)
\end{equation}
where  $\ln_+(x)=\max(\ln x,0)$ and $g$ is a bounded function over $D\times D$. In that case, one can show that $X$ lives in the space of distributions: this just means that for all smooth function $\varphi$ with compact support the integral $\int_D \varphi(x) X(x) dx$ makes sense. In fact, even if we will not discuss this here, GMC measures associated to log-correlated $X$, namely with covariance \eqref{logcov},  appear in many other fields among which: mathematical finance (see \cite{cf:BaKoMu} for a review), 3d turbulence \cite{ChRoVa}, decaying Burgers turbulence \cite{FLR}, the extremes of log-correlated Gaussian fields \cite{bramson,biskup,madaule}, the glassy phase of disordered systems \cite{CarDou,Fyo,rosso,MRV}  or the eigenvalues of Haar distributed random matrices \cite{webb}. However, we will focus in these notes on applications to LQG.

The purpose of these Les Houches lecture notes is twofold: first,  give a rigorous definition of measures of the type \eqref{firstdef} and review some of their main properties. Emphasis will be put on explaining the main ideas and not on giving rigorous proofs. Second, we will show how to use these measures to define Polyakov's 1981 theory of LQG \cite{Pol} on the Riemann sphere; in this specific case, one can identify LQG with Liouville quantum field theory (LQFT, hereafter). Here, emphasis will be put on explaining the construction of the so-called Liouville measures and explaining their (conjectured) relation with random planar maps: the construction will rely on the previous section on GMC.

\subsubsection*{Notations}

We will denote by $|.|$ the standard Euclidean metric, i.e. $|x-y|$ will denote the distance between two points $x$ and $y$. Also, if $A$ is some set then $|A|$ will stand for the Euclidean volume of $A$. It should be clear from the context whic convention is used for $|.|$. The Eucliden ball of center $x$ and radius $r>0$ will be denoted $B(x,r)$. The standard Lebesgue measure will be $dx$ in section 2; however, in section 3 on LQG, we will work exclusively in $2d$ so the Lebesgue measure will be denoted $dz$.

In these lecture notes, we will only study the theory of GMC in the case where $D$ is some subdomain of $\R^d$, $\sigma$ is  a Radon measure of the form $f(x) dx$ with $dx$ the Lebesgue measure, $f$ a nonnegative $L^1(dx)$ function and $X$ has a covariance kernel of log-type \eqref{logcov}. The underlying probability space will be $(\Omega, \mathcal{F}, \P)$ and we will denote $\E[.]$ the associated expectation. The vector space of $p$ integrable random variables with $p \geq 1$ will be denoted $L^p$. We will call a function $\theta:\R \to \R$ a smooth mollifier if $\theta$ is $C^\infty$ with compact support and such that $\int_{\R^d}  \theta(x)dx=1$. We will use $\theta$ to regularize the field $X$ by convolution; we will denote by $f \ast g$ the convolution between two distributions $f$ and $g$. When $\theta$ is smooth, the convolution $X\ast \theta$ is in fact $C^\infty$ and in particular the exponential of $X\ast \theta$ is well defined.    

In section 2, we will also consider centered Gaussian fields $Y,Z$ with  continuous covariances kernels and which are almost surely continuous.

\subsubsection*{Acknowledgments}

We would like to thank D. Chelkak for useful discussions on the Ising model and C. Hongler, F. David for the images. We also thank Y. Huang for reading carefully a prior draft of these lecture notes.

\section{Gaussian multiplicative chaos}

Before explaining the construction of the GMC measures, we first give a few reminders  on Gaussian vectors and processes.

\subsection{Reminder on Gaussian vectors and processes}

Here, we recall basic properties of Gaussian vectors and processes that we will need in these lecture notes. The first one is the Girsanov transform:

 \begin{theorem}{\bf Girsanov theorem }\label{th:Girsanov}

Let $(Y(x))_{x \in D}$ be a smooth centered Gaussian field with covariance kernel $K$ and $Y$ some Gaussian variable which belongs to the $L^2$ closure of the subspace spanned by $(Y(x))_{x \in D}$. Let $F$ be some bounded function defined on the space of continuous functions. Then we have the following identity
\begin{equation*}
\E[  e^{Y-\frac{\E[Y^2]}{2}}   F( (Y(x))_x   )   ]= \E[    F( (Y(x)   +E[Y Y(x)])_x   )   ]
\end{equation*}
\end{theorem}  
Though we state the Girsanov theorem under the above form, it is usually stated in the following equivalent form: under the new probability measure $ e^{Y-\frac{\E[Y^2]}{2}} d \P $, the field $(Y(x))_{x \in D}$ has same law as the (shifted) field $  (Y(x)   +E[Y Y(x)])_{x \in D}$ under $\P$.

We will also need the following beautiful comparison principle first discovered by Kahane:

 \begin{theorem}{\bf Convexity inequalities. [Kahane, 1985]. }\label{th:compar}

Let $(Y(x))_{x \in D}$ and $(Z(x))_{x \in D}$ be continuous centered Gaussian fields such that
\begin{equation*}
\E[ Y(x) Y(y)   ] \leq \E[ Z(x) Z(y)   ]. 
\end{equation*}
Then for   all convex (resp. concave) functions $F:\R_+\to \R$ with at most polynomial growth at infinity
\begin{equation}\label{eq:compar}
\E\Big[F \left(    \int_D e^{   Y(x)- \frac{\E[Y(x)^2]}{2}   } \sigma(dx)   \right ) \Big]\leq (\text{resp. } \geq)\, \E\Big[F \left(    \int_D e^{   Z(x)- \frac{\E[Z(x)^2]}{2}   } \sigma(dx)   \right ) \Big].
\end{equation}
\end{theorem}

\subsection{Construction of the GMC measures}

In this section, we will state a quite general theorem which will be used as definition of the GMC measure. The idea to construct a GMC measure is rather simple and standard: one defines the measure as the limit as $\epsilon$ goes to $0$ of $ c_\epsilon e^{  \gamma X_\epsilon} \sigma(dx)$ where $X_\epsilon$  is a sequence which converges to $X$ as $\epsilon$ goes to $0$ and $c_\epsilon$ is some normalization sequence which ensures that the limit is non trivial.
\begin{theorem}\label{th:existencechaos}
Let $\theta$ be a smooth mollifier. Set $X_\epsilon=X \ast \theta_\epsilon$ where $X$ has a covariance kernel of log-type \eqref{logcov} and $\theta_\epsilon= \frac{1}{\epsilon^d}  \theta ( \frac{.}{\epsilon} )$. The random measures
\begin{equation*}
M_{\epsilon,\gamma } (dx)=    e^{  \gamma X_\epsilon- \frac{ \gamma^2 \E[  X_\epsilon(x)^2  ]}{2} } \sigma(dx)
\end{equation*}
converge in probability in the space of Radon measures (equipped with the topology of weak convergence) towards a random measure $M_\gamma$. The random measure does not depend on the mollifier $\theta$. If $\sigma(dx)=f(x)dx$ with $f>0$, the measure $M_\gamma$ is different from $0$ if and only if $\gamma<\sqrt{2d}$.
\end{theorem}

\proof
For simplicity, we will prove the above theorem in the simple case where $\gamma<\sqrt{d}$, the so-called $L^2$ case. It is no restriction to suppose $f=1$ in the proof (the proof works the same with general $f$). Let $\theta$ be some smooth mollifier and $X_\epsilon=X \ast \theta_\epsilon$. For all compact $A$, we have by Fubini
\begin{equation*}
\E[  M_{\epsilon,\gamma } (A)   ]   =     \int_A  \E[ e^{  \gamma X_\epsilon(x)- \frac{ \gamma^2 \E[  X_\epsilon(x)^2  ]}{2} }  ]   dx= |A|.  
\end{equation*}
Hence, we see that the average of   $M_{\epsilon,\gamma } (A)$ is constant and equal to the Lebesgue volume of $A$: this explains the normalization term $\frac{ \gamma^2 \E[  X_\epsilon(x)^2  ]}{2}$ in the exponential. By a simple computation, one can show that for all $\epsilon' \leq \epsilon$ there exists global constants $c,C>0$ such that 
\begin{equation}\label{fundinequality}
c+ \ln \frac{1}{|y-x|+\epsilon}  \leq \E[  X_{\epsilon'}(x) X_\epsilon(y)   ]  \leq C + \ln \frac{1}{|y-x|+\epsilon}
\end{equation} 
One can notice that the bounds in the above inequality are independent of the smaller scale $\epsilon'$. Hence, using Fubini, we get that for all compact $A$
\begin{align*}
\E[  M_{\epsilon,\gamma } (A) ^2  ]   &  = \E  \left [   \left ( \int_A e^{  \gamma X_\epsilon(x)- \frac{ \gamma^2 \E[  X_\epsilon(x)^2  ]}{2} } dx \right )^2  \right ] \\
& = \int_A \int_A  \E \left [   e^{  \gamma (X_\epsilon(x)+X_\epsilon(y))- \frac{ \gamma^2 \E[  X_\epsilon(x)^2  ]}{2}-\frac{ \gamma^2 \E[  X_\epsilon(y)^2  ]}{2} }  \right  ]  dx dy  \\
& =  \int_A \int_A  e^{\gamma^2  \E[  X_\epsilon(x) X_\epsilon(y)   ] } dx dy   \\  
& \underset{\epsilon \to 0}{\rightarrow}  \int_A \int_A  e^{\gamma^2  K(x,y)   } dx dy,
\end{align*} 
 where the last convergence is a consequence of the simple convergence of $\E[  X_\epsilon(x) X_\epsilon(y)   ] $ towards $K$ for $x \not = y$ and the dominated convergence theorem using \eqref{fundinequality} (the condition $\gamma^2<d$ ensures the integrability of $e^{\gamma^2 K(x,y)}$).

Now, along the same lines (using Fubini), one can expand for $\epsilon'< \epsilon$ the quantity $\E[    (M_{\epsilon,\gamma } (A) -M_{\epsilon',\gamma } (A) )^2   ]$ and show that
\begin{align*}
& \E[  (M_{\epsilon,\gamma } (A) -M_{\epsilon',\gamma } (A) )^2  ]  \\
 &  =  \int_A \int_A  e^{\gamma^2  \E[  X_\epsilon(x) X_\epsilon(y)   ] } dx dy+  \int_A \int_A  e^{\gamma^2  \E[  X_{\epsilon'}(x) X_{\epsilon'}(y)   ] } dx dy - 2  \int_A \int_A  e^{\gamma^2  \E[  X_{\epsilon'}(x) X_\epsilon(y)   ] } dx dy    \\  
& \underset{\epsilon', \epsilon \to 0}{\rightarrow}  \int_A \int_A  e^{\gamma^2  K(x,y)   } dx dy+ \int_A \int_A  e^{\gamma^2  K(x,y)   } dx dy-2 \int_A \int_A  e^{\gamma^2  K(x,y)   } dx dy  \\
& =0  \\
\end{align*} 
hence $ ( M_{\epsilon,\gamma } (A))_{\epsilon>0} $ is a Cauchy sequence.  

Let  $\bar{\theta}$ be another smooth mollifier and let $\bar{M}_{\epsilon,\gamma } (dx)=    e^{  \gamma \bar{X}_\epsilon- \frac{ \gamma^2 \E[  \bar{X}_\epsilon(x)^2  ]}{2} } dx$ with $\bar{X}_\epsilon=X \ast \bar{\theta}_\epsilon$. Along the same lines as previously, one can show that $M_{\epsilon,\gamma } (A) -\bar{M}_{\epsilon,\gamma } (A)$ converges to $0$. 

In conclusion, we have shown that for all compact $A$, the variable $M_{\epsilon,\gamma } (A)$ converges in $L^2$ to some random variable $Z(A)$ of mean $|A|$,  and the limit $Z(A)$ does not depend on the smooth mollifier $\theta$. Using standard results of the theory of random measures (see \cite{daley}), one can show that there exists a random measure version $M_\gamma$ of the variables $Z(A)$ such that in fact $M_{\epsilon,\gamma }$ converges in probability in the space of random measures (equipped with the weak topology) towards $M_\gamma$. Of course, $M_\gamma$ is non trivial since for all compact $A$ we have $\E[M_\gamma(A)  ]=|A|$.

Now for the case $\gamma \in [\sqrt{d}, \sqrt{2d}[$, the above $L^2$ computations no longer converge and one must use more refined techniques to show convergence: we refer to Berestycki's approach \cite{Be} for a simple proof in that case.

\qed

\begin{figure}[h]
\centering
\includegraphics[width=0.8\linewidth]{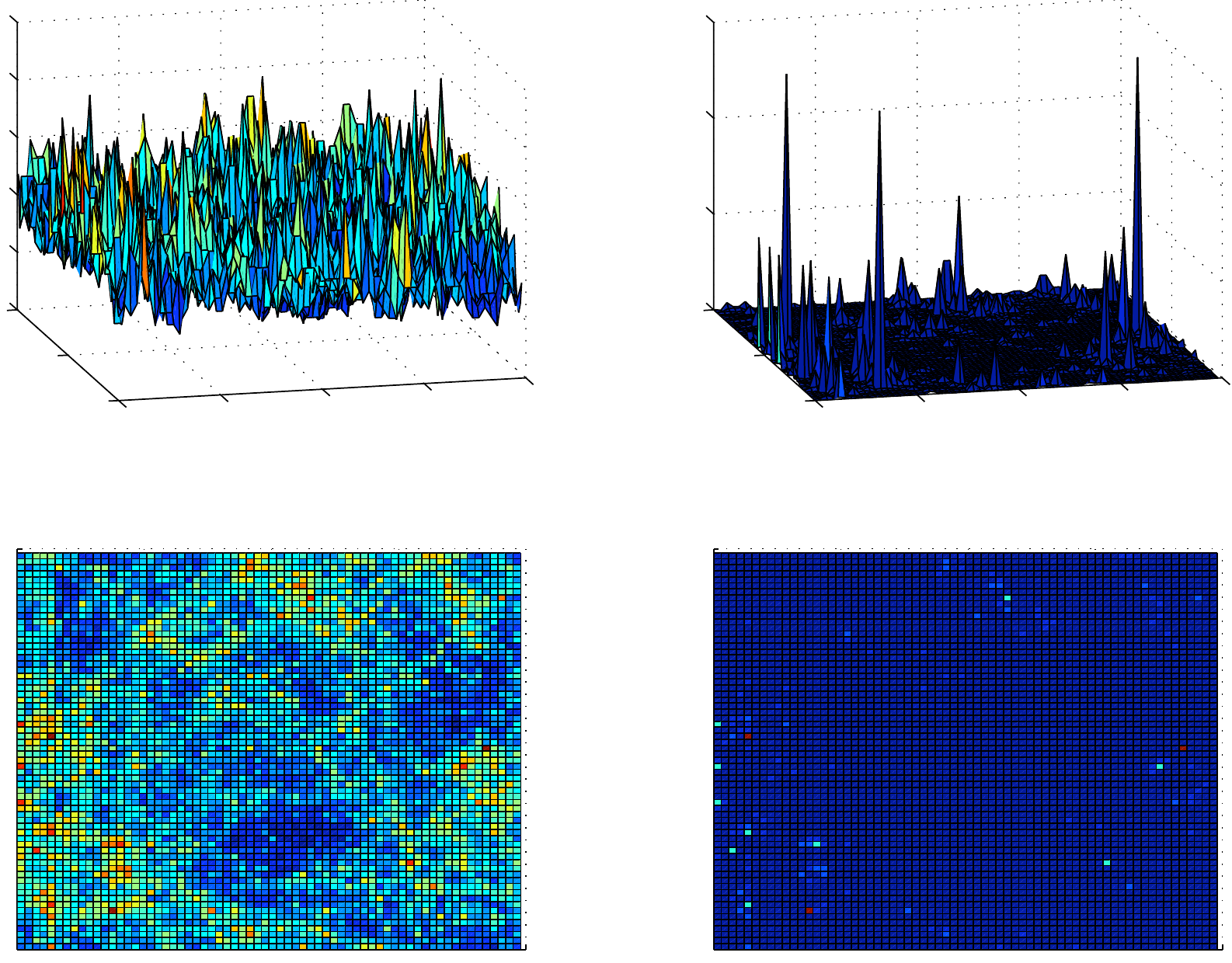}
\caption{Two examples of GMC measures.  Left: weak parameter $\gamma$. Right: $\gamma$ close to $2$.}\label{figchaos}
\end{figure}

\subsubsection{A brief historic on the construction of the GMC measures}
In fact, the above convergence result could be strengthened to more general cut-off approximations $X_\epsilon$ of the field $X$. However, for the sake of simplicity, we only stated the theorem with approximations $X_\epsilon$ of the form $X \ast \theta_\epsilon$. Before stating important properties of the measures $M_\gamma$, let us briefly review the historics of the above theorem. In his 1985 founding paper, Kahane defined the GMC measures by using a sequence of discrete approximations $X_n$ to $X$: he considered the simplified assumption that the random functions $(X_{n+1}(\cdot)-X_n(\cdot))_n$ are  independent\footnote{Kahane's motivation was the rigorous construction of Mandelbrot's limit lognormal model in turbulence defined in \cite{cf:Man}. Part of Mandelbrot's work \cite{cf:Man} is rigorous; Hoegh-Krohn also proved in \cite{Hoeg} similar results to \cite{cf:Man} around the same time.} . Within this framework, he defined the GMC measure as the almost sure limit of $M_{n,\gamma } (dx)=    e^{  \gamma X_n(x)- \frac{ \gamma^2 \E[  X_n(x)^2  ]}{2} } \sigma(dx)$ and showed that the law of the limiting measure does not depend on the sequence $X_n$. Around 20 years later, Robert-Vargas \cite{cf:RoVa} proved a weak form of theorem \ref{th:existencechaos} by showing convergence in law of $M_{\epsilon,\gamma } (dx)$. Duplantier-Sheffield \cite{cf:DuSh} proved theorem \ref{th:existencechaos} in the special case where $X$ is the GFF and $X_\epsilon$ is a circle average\footnote{Duplantier-Sheffield call this specific GMC measure the Liouville measure; in these notes, we choose a different convention for the terminology Liouville measure.} (this work was followed by the work of Chen-Jakobson \cite{CJ} where the authors adapted arguments fom \cite{cf:DuSh} to the $4d$ case). Recently, the convergence in law proved in \cite{cf:RoVa}  was reinforced to a convergence in probability by Shamov \cite{shamov}; the work of Shamov \cite{shamov}, which relies on abstract Gaussian space theory, is in fact quite general and does not concern just log-correlated $X$. Finally, let us mention that other works have now also established theorem \ref{th:existencechaos} by rather elementary methods: see Berestycki \cite{Be} and Junnila-Saksman \cite{JuSak} (this work is also interesting because it extends the theory to the critical case $\gamma=\sqrt{2d}$ where one can define a modified GMC theory; however, we will not consider  the critical case $\gamma=\sqrt{2d}$ in these notes). Berestycki's work \cite{Be} is probably a very good starting point for someone who wants to learn GMC theory.

\subsection{Main properties of the GMC measures}

Now, we turn to some important properties of the GMC measures which we will need in our study of LQG.

\subsubsection{Existence of moments and multifractality}

\begin{theorem}\label{th:existencemoments}
For $\gamma< \sqrt{2d}$, let $M_\gamma$ be a GMC measure associated to a log-correlated field $X$ with covariance \eqref{logcov} and $\sigma(dx)=f(x)dx$ with bounded $f$. Then, for $O \subset D$ an open ball we have
\begin{equation*}
\E[ M_\gamma (O)^q  ] < \infty
\end{equation*}
  if and only if $q \in ]- \infty, \frac{2d}{\gamma^2}[$.
\end{theorem}

We will not prove this theorem here: we refer to \cite{cf:RoVa} for a proof. Now, we turn to the multifractal scaling of the measure. This is the content of:
\begin{proposition}
For $\gamma< \sqrt{2d}$, let $M_\gamma$ be a GMC measure associated to a log-correlated field $X$ with covariance \eqref{logcov} and $\sigma(dx)=f(x)dx$ with bounded continuous $f$. Then for all $x$ and all $q \in ]- \infty, \frac{2d}{\gamma^2}[$, there exists some constant $C_x>0$ (which depends also on $f$, $q$ and the exact form of the kernel $K$ in  \eqref{logcov}) such that
\begin{equation}\label{scalinglocal}
\E[   M_\gamma (B(x,r))^q  ]  \underset{r \to 0}{\sim} C_x r^{\zeta(q)}
\end{equation}
where $\zeta(q)= (d+\frac{\gamma^2}{2})q- \frac{\gamma^2 q^2}{2}$ is called the structure function of $M_\gamma$. 
\end{proposition}

The above proposition implies that the GMC measure associated to a log-correlated field $X$ exhibits multifractal behaviour, i.e. the measure is not scale invariant but rather is locally H\"older around each point. The H\"older exponent depends on the point (for more on the so-called multifractal formalism, see the next subsection).   More generally, one can take as a definition that a random measure  satisfying \eqref{scalinglocal} where $\zeta$ is a strictly concave function is a multifractal measure.

\subsubsection{Multifractal formalism}

Now, we turn to the multifractal formalism of the measures $M_\gamma$. The measures $M_\gamma$ are multifractal in the sense that the regularity of the measure around a point $x \in D$ depends on the point $x$: this can easily be seen on figure \ref{figchaos}. Multifractal formalism is a general theory to study the regularity of measures like $M_\gamma$ around each point:  for more background on this see section 4 in \cite{review}.

For $\gamma^2<2d$ and $q \in ]0, \frac{\sqrt{2d}}{\gamma}[$, we consider the following set:
 \begin{equation*}
 K_{\gamma,q}= \left\{ x \in D; \: \underset{\epsilon \to 0}{\lim} \: \frac{\ln M_{\gamma}(B(x,\epsilon))}{\ln \epsilon} = d+(\frac{1}{2}-q)\gamma^2 \right\}.
 \end{equation*}
In words, the set $K_{\gamma,q}$ is made of the points $x$ such that $M_\gamma (B(x,r))  \underset{r \to 0}{\approx} r^{d+(\frac{1}{2}-q)\gamma^2}$.
We can state the following theorem:
\begin{theorem}
The set $K_{\gamma,q}$ has Hausdorff dimension $d-\frac{\gamma^2 q^2}{2}$.
\end{theorem}

 In fact, the same theorem holds with the set $\bar{K}_{\gamma,q}$ defined by
 \begin{equation*}
 \bar{K}_{\gamma,q}=  \left\{ x \in D; \: \underset{\epsilon \to 0}{\lim} \: \frac{ X_\epsilon(x)}{-\ln \epsilon} = \gamma q \right\},
 \end{equation*}
where $X_\epsilon =X \ast \theta_\epsilon$ with $\theta$ any smooth mollifier. The reason is that it is useful to have in mind the following approximation 
\begin{equation}\label{approxball}
M_\gamma (B(x,r)  )  \underset{r \to 0}{\approx}  r^d  e^{\gamma X_r(x)-\frac{  \gamma^2  \E[ X_r(x)^2 ]  }{2} }
\end{equation}
where here $a_r \approx b_r$ means that the ratio $a/b$ is a (random) constant $C_r$ of order $1$, i.e. $\E[ C_r ]$ belongs to an interval $[c,C]$ with $c,C>0$ independent of $r$. The main difficulty in our context is that the random constant $C_{x,r}$ for the ratio of both sides in \eqref{approxball} really also depends on $x$ so the above approximation can not be used directly but is rather a guideline to get intuition on the behaviour of $M_\gamma$. In our case, if we assume $C_{x,r}=1$, then the sets $K_{\gamma,q}=\bar{K}_{\gamma,q}$ are the same (however we stress that rigorously these two sets are not the same). Finally, following the terminology of Hu-Miller-Peres \cite{HMP}, a point $x$ which belongs to $\bar{K}_{\gamma,1}$ is nowadays called a $\gamma$-thick point.

Now, among the sets $K_{\gamma,q}$ (and $\bar{K}_{\gamma,q}$), the set $K_{\gamma,1}$ (resp. $\bar{K}_{\gamma,1}$) is of particular importance for $M_\gamma$ since it is the set on which the measure $M_\gamma$ "lives". More specifically, we have
\begin{equation}\label{strongthick}
M_\gamma  (^{c} K_{\gamma,1} \cup {}^{c} \bar{K}_{\gamma,1} )=0
\end{equation}
In the modern terminology of \cite{HMP}, one says that $M_\gamma$ lives on the $\gamma$-thick points of $X$. This property was proved by Kahane in his  seminal paper \cite{cf:Kah}. Here, we will show a slightly weaker result, namely that:
\begin{equation}\label{weakthick}
M_\gamma \left ( ^c \left\{ x \in D; \: \underset{n \to \infty}{\lim} \: \frac{ X_{\frac{1}{2^n}}(x)}{n \ln2 } = \gamma  \right\} \right )=0
\end{equation} 
The only difference with $\bar{K}_{\gamma,1}$ is that we restrict the limit in $\bar{K}_{\gamma,1}$ to a dyadic sequence (in fact, with little effort, one can reinforce \eqref{weakthick} to prove \eqref{strongthick}). 

\emph{Proof of \eqref{weakthick}}:

We introduce $\eta>0$ and a compact set $A$. We have for all $n \leq p$ and by using the Girsanov theorem \ref{th:Girsanov} that
\begin{align*}
& \E[\int_A  \ind_{ \left\{ x \in D; \:   \frac{X_{2^{-n}}(x)}{n \ln 2} \in {}^c [\gamma-\eta, \gamma+\eta]  \right\}}    e^{  \gamma X_{2^{-p}}(x)  - \frac{\gamma^2}{2} \E[  X_{2^{-p}}(x)^2  ]  } dx ]   \\
& = \int_A  \E[  \ind_{ \left\{ x \in D; \:   \frac{X_{2^{-n}}(x)}{n \ln 2} \in {}^c [\gamma-\eta, \gamma+\eta]  \right\}}    e^{  \gamma X_{2^{-p}}(x)  - \frac{\gamma^2}{2} \E[  X_{2^{-p}}(x)^2  ]  }  ]  dx  \\
& = \int_A  \E[  \ind_{ \left\{ x \in D; \:   \frac{X_{2^{-n}}(x)+  \gamma \E[ X_{2^{-n}}(x) X_{2^{-p}}(x)   ]}{n \ln 2} \in {}^c [\gamma-\eta, \gamma+\eta]  \right\}}     ]  dx  \\
& \approx  \int_A \P   \left (     \frac{X_{2^{-n}}(x)   }{n \ln 2} \in {}^c [-\eta,  \eta]    \right )  dx
\end{align*} 
Now, since $ X_{2^{-n}}(x)$ is a Gaussian of variance roughly equal to $n \ln 2$ by \eqref{fundinequality}, we get that
\begin{equation*}
\P   \left (     \frac{X_{2^{-n}}(x)   }{n \ln 2} \in {}^c [-\eta,  \eta]  \right  ) \leq 2 e^{- n \eta^2 \frac{(\ln 2)^2}{2}}
\end{equation*}
 Therefore, by taking the limit $p \to \infty$ in the above considerations, we get that there exists $C>0$
\begin{equation*}
 M_\gamma \left (  \left\{ x \in D; \:  \frac{ X_{2^{-n}}(x)}{n \ln2 } \in {}^c [\gamma-\eta, \gamma+\eta]  \right\} \right ) \leq C e^{- n \eta^2 \frac{(\ln 2)^2}{2}}
\end{equation*}
One can easily deduce from this by a Borell-Cantelli type argument that 
\begin{equation*}
 M_\gamma \left (  \cap_N \cup_{n \geq N}  \left\{ x \in D; \:  \frac{ X_{2^{-n}}(x)}{n \ln2 } \in {}^c [\gamma-\eta, \gamma+\eta]  \right\} \right )=0
\end{equation*}
Since the result is valid for all $\eta>0$, we get \eqref{weakthick}.

\qed

\subsubsection{The first Seiberg bound}

In this subsection, we state and prove a theorem we will need to define LQG: indeed, we will see that it corresponds to the so-called Seiberg bound in LQG. We have the following 
\begin{lemma}\label{firstseib}
Let $\alpha \in \R$ and $x \in D$. We have
\begin{equation*}
\int_{B(x,1)}  \frac{1}{|y-x|^{\alpha \gamma}}  M_{\gamma} (dy) < \infty, \quad a.s.
\end{equation*}
if and only if $\alpha < \frac{d}{\gamma}+ \frac{\gamma}{2}$.
\end{lemma}

\proof

We only prove the if part; for the only if part, we refer to \cite{DKRV}. With no loss of generality, we suppose that $x=0$.  We consider $\eta \in ]0,1[$. We have 
\begin{align*}
\E   \left [ \left (  \int_{B(0,1)}  \frac{1}{|y|^{\alpha \gamma}}  M_{\gamma} (dy) \right )^\eta \right ] & \leq \sum_{n=1}^\infty \E \left [ \left ( \int_{\frac{1}{2^n}  \leq |y| \leq \frac{1}{2^{n-1}} }  \frac{1}{|y|^{\alpha \gamma}}  M_{\gamma} (dy) \right )^\eta\right ]  \\
& \leq \sum_{n=1}^\infty   2^{n \alpha \eta \gamma } \E[ (  M_{\gamma} (  \lbrace y;  \frac{1}{2^n}  \leq |y| \leq \frac{1}{2^{n-1}}  \rbrace) )^\eta ]  \\
& \leq \sum_{n=1}^\infty   2^{n \alpha \eta \gamma } \E[ (  M_{\gamma} (  \lbrace y;  |y| \leq \frac{1}{2^{n-1}}  \rbrace) )^\eta ]  \\
& \leq C \sum_{n=1}^\infty   2^{n \alpha  \gamma  \eta} 2^{-n \zeta (\eta)}, \\
\end{align*}   
where recall that $\zeta(q)= (d+\frac{\gamma^2}{2})q- \frac{\gamma^2 q^2}{2}$. Now, since $\alpha < \frac{d}{\gamma}+ \frac{\gamma}{2}$,  one can choose $\eta>0$ small such that $ \alpha \gamma \eta -\zeta (\eta)<0$ hence we get the conclusion.

\qed

\section{Liouville quantum gravity on the Riemann sphere}

Now, in the second part of these notes, we show how to use GMC theory to construct Liouville Quantum Gravity (LQG) on the Riemann sphere. LQG was introduced in Polyakov's seminal 1981 paper \cite{Pol}. In the paper \cite{Pol}, Polyakov builds a theory of summation of 2d-random surfaces in the spirit of Feynman's theory of summation of random paths. On the Riemann sphere, LQG is in fact equivalent to Liouville quantum field theory (LQFT); for a complete review on LQFT in the physics literature, we refer to Nakayama \cite{nakayama}. However, LQG is a general theory of random surfaces which can be defined on any 2d-surface. In the case of higher genus surfaces, LQFT is a building block of LQG and they are not equivalent.  For the sake of simplicity, we will restrict ourselves here to the case of the sphere where we identify LQG and LQFT: in this context, we explain the construction of LQFT following David-Kupiainen-Rhodes-Vargas \cite{DKRV}.

LQFT is not only a quantum field theory but since it has extra symmetries it is also a conformal field theory (CFT). Quantum field theory and conformal field theory is a very wide topic in mathematical physics which can be approached in different ways: by algebraic methods, geometric methods and probabilistic methods. Of course, all these approaches can be related but for reasons of simplicity (and the knowledge of the authors!) we will restrict to the probabilistic setting. Before we describe the theory, we first give a brief introduction to what is a CFT on the Riemann sphere. Then, we introduce a few notations and definitions from elementary Riemannian geometry.

\subsection{Elementary Riemannian geometry on the sphere}\label{sub:geom}
We consider the standard Riemann sphere $\mathbb{S}= \C \cup \lbrace \infty \rbrace$. The Riemann sphere $\mathbb{S}$ is just the complex plane $\C$ with a point at infinity and is obtained as the image of the standard $2d$ sphere by stereographic projection. We equip $\mathbb{S}$ with the standard round metric.  On $\mathbb{S}$, the round metric is given in Riemannian geometry notations  by $g(z) |dz|^2$ where $g(z)=\frac{4}{(1+|z|^2)^2} $. This means that the length $\mathcal {L}(\sigma)$ of a curve $\sigma: [0,1] \to \mathbb{S}$ is given by 
\begin{equation*}
\mathcal {L}(\sigma)= \int_0^1   g(\sigma(t))^{1/2} |\sigma'(t)| dt.
\end{equation*}  
One then gets the distance between two points $z_1,z_2 \in \mathbb{S}$ by taking the infimum of $\mathcal {L}(\sigma)$ over all curves $\sigma$ which join $z_1$ to $z_2$. The volume form is simply given by the measure $g(z) dz$ where $dz$ is the Lebesgue measure on $\R^2$ (by using polar coordinates, it is easy to see that $\int_{\mathbb{S}} g(z) dz=4 \pi$, thereby recovering the well known fact that the surface of the sphere is $4 \pi$!). In this context, one can of course do differential calculus and $C^k$ functions on $\mathbb{S}$ are just functions $\phi$ defined on $\C$ which are such that $\phi$ is $C^k$ on $\C$ and $z \mapsto  \phi(\frac{1}{z})$ admits a continuous extension on $\C$ which is $C^k$. The gradient $\nabla_g$ of a function $\phi$ is given by the simple formula 
\begin{equation*}
\nabla_g \phi(z) =\frac{1}{g(z)}\nabla_z \phi (z). 
\end{equation*}
where $\nabla_z$ is the standard Euclidean gradient on $\C$.
 Finally, the (Ricci) curvature $R_g$ is given by 
 \begin{equation*}
 R_g(z)= - \frac{1}{g(z)}\Delta_z \ln g(z),
 \end{equation*}
 where $\Delta_z$ is the standard Euclidean Laplacian. In the specific case of the round metric ($g(z)=\frac{4}{(1+|z|^2)^2} $), one finds by a simple computation a constant curvature $R_g=2$.

\subsection{An introduction to CFT on the Riemann sphere}\label{sub:introCFT}

The general formalism of CFT was built in the celebrated 1984 work of Belavin-Polyakov-Zamolodchikov \cite{BPZ}. Here we give an elementary (and incomplete) exposition of this formalism. A CFT on the Riemann sphere is usually defined by: 
\begin{enumerate}
\item
a real parameter $c_{CFT}$ called the central charge 
\item
 (primary) local fields $(\phi_\alpha)_{\alpha \in \mathcal{A}}$ defined in the complex plane $\C$.
 \item prescribed symmetries (conformal covariance, diffeomorphism invariance, Weyl anomaly: see Gawedzki's lecture notes \cite{gaw} for further details): see equality \eqref{eq:confcov} below for the conformal covariance statement.
 \end{enumerate}
 
  It is not obvious to give a simple definition of the central charge but we will see in the example of LQFT how it appears. For now, let us just mention that the central charge of a CFT determines the symmetries of the theory; however, it is very important to stress that two CFTs with same central charge can be very different because the set of primary local fields plays an essential role too. In a CFT theory, what makes sense are the correlation functions $<\phi_{\alpha_1}(z_1) \cdots \phi_{\alpha_n}(z_n)>$ at non coincident points $z_i$ (i.e. $z_i \not = z_j$ for $i \not = j$) and for certain values of the $\alpha_i$ where $<.>$ should be viewed as some underlying  measure (however, this is a view as the measure does not necessarily exist). The correlation functions of primary local fields satisfy the following conformal covariance: if $\psi$ is a M\"obius transform on the sphere $\mathbb{S}$, i.e. $\psi(z)= \frac{az+b}{cz+d}$ where $a,b,c,d \in \C$ are such that $ad-bc=1$, then
\begin{equation}\label{eq:confcov}
<\phi_{\alpha_1}(\psi(z_1)) \cdots \phi_{\alpha_n}(\psi(z_n))>= \prod_{i=1}^n |\psi'(z_i)|^{-2 \Delta_{\alpha_i}} <\phi_{\alpha_1}(z_1) \cdots \phi_{\alpha_n}(z_n)>
\end{equation}      
where the real number $\Delta_{\alpha_i}$ is called the conformal weight of the field $\phi_{\alpha_i}$.  One of the successes of CFT is that it describes (conjecturally in mathematical standards) the scaling limit of correlation functions of statistical physics models at critical temperature. It is a major program in mathematical physics to make these predictions from CFT rigorous mathematical statements.  

In some cases, one can also define $\phi_{\alpha}$ in a strong sense as a random distribution in the sense of Schwartz: in that case, if $\mathcal{D}$ denotes the set of smooth functions with compact support one can consider the random distribution $\varphi \in \mathcal{D} \to \int_{\C} \phi_\alpha(x) \varphi(x) dx  $. In this case, the underlying  measure really exists (this will be the case for some but not all primary local fields in the two examples we will consider in these notes: LQFT and the Ising model at critical temperature) and one can compute the moments of the variable $  \int \phi_\alpha(z) \varphi(z) dz$ (if they exist) in terms of the correlation functions by the following obvious formula    
\begin{equation*}
 < ( \int_{\C} \phi_{\alpha}(z) \varphi(z) dz )^n>= \int_{\C}  \cdots \int_{\C} <\phi_{\alpha}(z_1) \cdots \phi_{\alpha}(z_n)> \varphi(z_1) \cdots  \varphi(z_n) \: dz_1 \cdots dz_n.  
\end{equation*}      
Hence, in many cases, the correlation functions determine the joint laws of the collection $\left  ( \int_{\C} \phi_\alpha(z) \varphi(z) dz \right )_{\varphi \in  \mathcal{D}}$.

\subsection{Introduction to LQFT on the Riemann sphere}

LQFT is a family of CFTs parametrized by two constants $\gamma \in ]0,2]$  and $\mu >0$; in these notes, we will only consider the case $\gamma \in ]0,2[$. In the probabilistic setting, the goal of LQFT is to make sense of and compute as much as possible the following correlation functions which arise in theoretical physics under the following heuristic form:
\begin{equation*}
< e^{\alpha_1 X(z_1)} \cdots e^{\alpha_n X(z_n)  }  >: =\int e^{\alpha_1 X(z_1)} \cdots e^{\alpha_n X(z_n)} e^{-  S_L(X, g) } DX
\end{equation*}
where $DX$ is the "Lebesgue" measure on functions $\mathbb{S} \to \R$ and  $S_L$ is the Liouville action: 
\begin{equation}\label{Liouvilleaction}
S_L(X,g):= \frac{1}{4\pi} \int_{\S}\big(|\nabla_g X |^2(z)+QR_{g}(z) X(z)  +4\pi \mu e^{\gamma X(z)  }\big)\, g(z)dz 
\end{equation}
where recall that $g$ is the round metric, the constant $Q$ is defined by $Q=\frac{\gamma}{2}+\frac{2}{\gamma}$  and $\mu >0$. LQFT is therefore an {\bf interacting} quantum field theory where the interaction term is
\begin{equation}\label{interacterm}
\mu \int_{\S}  e^{\gamma X(z)  }\, g(z)dz.
\end{equation}
The positive parameter $\mu$, called the cosmological constant, is necessary for the existence of LQFT. However, a remarkable feature of LQFT is that the parameter $\gamma$ is the essential parameter of the theory as it completely determines the conformal properties of the theory (in CFT language, the parameter $\gamma$ determines the central charge: we will come back to this point later in more detail). Following the standard terminology of CFT (see previous chapter), the $e^{\alpha_i X(z_i)}$ are local primary fields (the conformal covariance property will be proved in the next chapter); in fact, in the context of LQFT, the  $e^{\alpha_i X(z_i)}$ are also called vertex operators.

It is a well known fact that the "Lebesgue measure" $DX$ does not exist since the space of functions $\mathbb{S} \to \R$ is infinite dimensional; however, it is a standard procedure in the probabilistic approach to quantum field theory (see Simon's reference book \cite{Simon} on the topic) to interpret the term $e^{- \frac{1}{4\pi} \int_{\R^2}  |\nabla_g X |^2(z) \, g(z)dz } DX$ as the Gaussian Free Field (GFF), i.e. the Gaussian field whose covariance is given by the Green function on $\S$. One way to see that this is the proper definition is to perform the following integration by parts
\begin{equation*}
e^{- \frac{1}{4\pi} \int_{\R^2}  |\nabla_g X |^2(z) \, g(z)dz } DX= e^{ \frac{1}{4\pi} \int_{\R^2}   X (z)  \Delta_g X(z) \, g(z)dz } DX
\end{equation*}
Formally, this corresponds to a Gaussian with covariance $2 \pi (- \Delta_g)^{-1}$.  In fact, thanks to the theory of probability, one can define the GFF rigorously with the following definition:

\begin{definition}
The GFF with vanishing mean on the sphere $X_g$ is the Gaussian field living in the space of distributions such that for all smooth functions $f,h$ on $\S$
\begin{equation*}
\E \left [    \left ( \int_{\S} f(z) X_g(z) g(z)dz \right )  \left (  \int_{\S} h(z') X_g(z') g(z')dz'  \right )  \right   ]   =   \int_{\S} \int_{\S} G_g(z,z') f(z) h(z')  g(z) g(z') dz dz' 
\end{equation*}
where $G$ is the Green function for the Laplacian on the sphere defined for all $z \in \S$ by
\begin{equation*}
-\Delta_g G(z,.)= 2 \pi (\delta_z-\frac{1}{4\pi}) , \; \; \int_{\S}  G_g(z,z' ) g(z') dz'=0 
\end{equation*} 
\end{definition}   
The random variable $X_g$ lives in the space of random distributions but in fact it exists in a Sobolev space and  $\int_{\S} f(z) X_g(z) g(z)dz $ makes sense for many functions $f$ (with less regularity than $C^\infty$). In particular, $ \int_{\S} X_g(z) g(z)dz$ makes sense and is equal to $0$ actually: this is why we call $X_g$ the GFF with vanishing mean on the sphere. It turns out that the Green function has the following explicit form on $\S$
\begin{equation*}
G_g(z,z')=  \ln  \frac{(1+|z|^2)^{1/2}(1+|z'|^2)^{1/2} }{|z-z'|},  
\end{equation*}
where recall that $|.|$ is the standard Euclidean distance. 

Now, in the spirit of probabilistic quantum field theory, since formally we have 
\begin{equation*}
e^{-  S_L(X, g) } DX= e^{ -\frac{1}{4\pi} \int_{\S}\big(QR_{g}(z) X(z)  +4\pi \mu e^{\gamma X(z)  }\big)\, g(z)dz  }    \times  e^{- \frac{1}{4\pi} \int_{\R^2}  |\nabla_g X |^2(z) \, g(z)dz } DX
\end{equation*}
and since we interpret $ e^{- \frac{1}{4\pi} \int_{\R^2}  |\nabla_g X |^2(z) \, g(z)dz } DX$ as the GFF measure, it is natural to interpret the formal measure $e^{-  S_L(X, g) } DX$ as follows, for all functions $F$ (up to a global constant) 
\begin{equation}\label{wrongdef}
\int F(X)e^{-  S_L(X, g) } DX=  \underset{\epsilon \to 0 }{ \lim}  \:   \E  [  F(X_g) e^{ - \frac{Q}{4 \pi} \int_{\mathbb{S}}  R_g(x) X_g(z) g(z) dz  -\mu \epsilon^{\gamma^2/2}  \int_{\mathbb{S}} e^{\gamma \bar{X}_{\epsilon,g}(z)} g(z) dz  }    ]  
\end{equation}
where  $\bar{X}_{\epsilon,g}$ is the average of $X_g$ in a ball of radius $\epsilon$ with respect to the metric $g$. However, there is something wrong with definition \eqref{wrongdef}; though it can be used to define a standard quantum field theory in the spirit of \cite{Simon}, it will lack symmetry to define a CFT. The reason is that we have not taken into account the contribution of constant functions in the Gaussian measure. This omission reflects in the fact that there is something arbitrary in the choice of $X_g$: indeed, $X_g$ has vanishing mean on the sphere but we could have chosen an other GFF on the sphere. In particular, $X_g$ is not conformally invariant since for all M\"obius transform $\psi$ the following equality holds in distribution:
\begin{equation*}
X_g \circ \psi - \int_{\mathbb{S}} (X_g \circ \psi (z)) g(z) dz= X_g.
\end{equation*}  
Now, the average $\int_{\mathbb{S}} X_g \circ \psi (z) g(z) dz$ is a Gaussian random variable which is non zero (unless $\psi$ is an isometry of $\mathbb{S}$) and hence  $X_g \circ \psi$ does not have the same distribution as $X_g$. One very natural way to get rid of this average dependence is to replace $X_g$ by $X_g+c$ where $c$ is distributed according to the Lebesgue measure (and stands for the mean value of the field). This leads to the following correct definition (up to some global constant)
\begin{equation}\label{rightdef}
\int F(X) e^{-  S_L(X, g) } DX=  \underset{\epsilon \to 0 }{ \lim}   \int_\R \E  [  F(X_g+c) e^{ - \frac{Q}{4 \pi} \int_{\mathbb{S}}  R_g(z)(X_g(z)+c) g(z) dz  -\mu \epsilon^{\gamma^2/2} e^{\gamma c}  \int_{\mathbb{S}} e^{\gamma \bar{X}_{\epsilon,g}(z)} g(z) dz  }    ]   dc
\end{equation}
A standard computations shows that 
\begin{equation}\label{devvar}
\E[   \bar{X}_{\epsilon,g}(z)  ^2 ] = \ln \frac{1}{\epsilon}- \frac{1}{2} \ln g(z)+ C+o(1)
\end{equation} 
where $C$ is some global constant, therefore the measure $\epsilon^{\gamma^2/2}  e^{\gamma \bar{X}_{\epsilon,g}(z)} g(z) dz$ converges to $e^{\frac{\gamma^2}{2} C}$ times the GMC measure $M_\gamma$ associated to $X_g$ and $g(z)dz$ which we write
\begin{equation}\label{GMCLQFT}
 M_\gamma(dz)=  e^{\gamma X_g(z)} g(z) dz.
\end{equation}
By the previous results on GMC theory, this GMC measure is well defined and non trivial. In the sequel, we will exclusively work with this GMC measure.

\subsection{Construction of LQFT}
  
 With the preliminary remarks of the previous subsection, we are ready to introduce the correlation functions of LQFT on the sphere and recover many known properties in the physics literature. In fact, it is standard in the physics literature to express the correlations of LQFT in the complex plane and therefore to shift the metric dependence of the theory in the field $X_g+c$: this simplifies many computations. Let us describe how to do so. If $\epsilon$ is small then a ball  $B_g(z,\epsilon)$ of centre $z$ and radius $\epsilon$ in the round metric $g$ is to first order in $\epsilon$ the same as an Euclidean ball $B(z, \frac{\epsilon}{g(z)^{1/2}})$ of centre $z$ and radius $\frac{\epsilon}{g(z)^{1/2}}$. Hence, the average $\bar{X}_{\epsilon,g}(z)$ (with respect to balls in the round metric) is roughly the same as $X_{\frac{\epsilon}{g(z)^{1/2}},g}(z)$ where $X_{\epsilon,g}(z)$ is the average of $X_g$ on an Euclidean ball of radius $\epsilon$. Finally, notice that we can write for all $\epsilon'>0$
 \begin{equation}\label{changemetric}
 (\epsilon')^{\gamma^2/2}  \int_{\mathbb{S}} e^{\gamma \bar{X}_{\epsilon',g}(z)} g(z) dz=    \int_{\mathbb{S}} \left ( \frac{ \epsilon'}{g(z)^{1/2}} \right )^{\gamma^2/2} e^{\gamma (\bar{X}_{\epsilon',g}(z)   + \frac{Q}{2} \ln g(z))}  dz 
 \end{equation} 
  where recall that $Q=\frac{\gamma}{2}+\frac{2}{\gamma}$. Since $\bar{X}_{\epsilon',g}(z) \approx X_{\frac{\epsilon'}{g(z)^{1/2}},g}(z)$, by making the change of variable $\epsilon=  \frac{ \epsilon'}{g(z)^{1/2}}$ in \eqref{changemetric}, it is not suprising that one can prove that the random measures   
\begin{equation*}
e^{\gamma (X_{\epsilon,g}(z)   + \frac{Q}{2} \ln g(z))}  dz 
\end{equation*}
  converge in probability as $\epsilon$ goes to $0$ towards $e^{\frac{\gamma^2}{2} C} M_\gamma(dz)$ where $C$ is defined by \eqref{devvar} and $M_\gamma$ is defined by \eqref{GMCLQFT}. Therefore, instead of working with $X_g+c$, we will work with the shifted field $\phi(z)= X_g(z)+c+ \frac{Q}{2} \ln g(z)$ and the approximations $\phi_\epsilon (z)= X_{\epsilon,g}+c+ \frac{Q}{2} \ln g(z)$. The field $\phi$ under the probability measure \eqref{rightdef} is called the {\bf Liouville field}. Finally, we set formally $V_\alpha(z)= e^{\alpha \phi(z)}$ and define the associated approximate  vertex operators
 \begin{equation}\label{approxvertex}
 V_{\alpha,\epsilon}(z):= \epsilon^{\alpha^2/2} e^{\alpha \phi_\epsilon (z) }.
 \end{equation}
The correlation functions of LQFT are now defined by the following formula
  \begin{equation}\label{defcorrels}
 < \prod_{i=1}^n V_{\alpha_i}(z_i) >:=  Z_{\textrm{GFF}}(g)  \:   \underset{\epsilon \to 0 }{ \lim}   \int_\R \E  [  \prod_{i=1}^n   V_{\alpha_i,\epsilon}(z_i)  e^{ - \frac{Q}{4 \pi} \int_{\mathbb{S}}  R_g(z)(X_g(z)+c)g(z) dz  -\mu \epsilon^{\gamma^2/2}  \int_{\mathbb{S}} e^{\gamma \phi_{\epsilon}(z)}  dz  }    ]   dc,
  \end{equation}  
  where one can notice the presence of the partition function of the GFF $ Z_{\textrm{GFF}}(g)$ given by  $\text{Det}\:  \Delta_g^{-1/2}$ where $\text{Det} \: \Delta_g$ is the standard determinant of the Laplacian (this determinant is in fact non trivial to define since the Laplacian is defined on an infinite dimensional space: see \cite{dubedat} for background). The constant $Z_{\textrm{GFF}}(g)$ is a global constant and plays no role here so it is not important to understand exactly how it is defined. For the readers who are unfamiliar with $\text{Det}\: \Delta_g$ they can take out this term in definition \eqref{defcorrels} and remember that it only plays a role in the Weyl anomaly formula (see proposition \ref{Weylanom} below).   
  
  Of course, it is crucial to enquire when the limit \eqref{defcorrels} exists. This is the object of the following:
  
 \begin{proposition}[\cite{DKRV}]\label{seibergprop} 
 The correlation functions \eqref{defcorrels} exist  and are not equal to $0$ if and only if the following Seiberg bounds hold
 \begin{equation}\label{Seiberg}
  \forall i, \; \alpha_i<Q \;  \; \;  \; \text{and}     \; \; \; \; \sum_{i=1}^n \alpha_i >2Q.
 \end{equation}
 In particular, the number of vertex operators $n$ must be greater or equal to $3$ for the correlation functions to exist and be non trivial.  
  If the Seiberg bounds hold  then we get the following expression (up to some multiplicative constant which plays no role and depends on the $C$ of \eqref{devvar}, $\alpha_1, \cdots, \alpha_n$ and $\gamma$)
  \begin{equation}\label{expcorrels}
  < \prod_{i=1}^n V_{\alpha_i}(z_i) >= Z_{\mathrm{GFF}}(g) \: e^{\frac{1}{2}\sum_{i \not = j} \alpha_i \alpha_j G_g(z_i,z_j)  } \:  \prod_{i=1}^n g(z_i)^{\frac{\alpha_i Q}{2}-\frac{\alpha_i^2}{4}} \: \Gamma (\frac{\sum_i \alpha_i-2Q}{\gamma} ,\mu) \:   \E  [  (  Z_{(z_i,\alpha_i)}(\S)^{-\frac{\sum_i \alpha_i-2Q}{\gamma} }   ] 
  \end{equation} 
 where $\Gamma (\frac{\sum_i \alpha_i-2Q}{\gamma} ,\mu)= \int_0^\infty u^{\frac{\sum_i \alpha_i-2Q}{\gamma}-1} e^{-\mu u} du$ and 
  \begin{equation*}
  Z_{(z_i,\alpha_i)}(dz)= e^{\gamma \sum_{i=1}^n \alpha_i G_{g}(z_i,z)} M_{\gamma} (dz).
  \end{equation*} 
   \end{proposition}
   
   \proof
   Here, we give a sketch of the proof of the if part of proposition \ref{seibergprop}: therefore, we suppose that the $(\alpha_i)_i$ satisfy the Seiberg bounds \eqref{Seiberg}. We denote $ < \prod_{i=1}^n V_{\alpha_i,\epsilon}(z_i) >$ the right hand side of  \eqref{defcorrels}. Since $R_g=2$ and $X_g$ has vanishing mean on the sphere one has   
  \begin{align*}
    < \prod_{i=1}^n V_{\alpha_i,\epsilon}(z_i) > /    Z_{\mathrm{GFF}}(g)& =     \int_\R \E  [  \prod_{i=1}^n   V_{\alpha_i,\epsilon}(z_i)  e^{ - \frac{Q}{4 \pi} \int_{\mathbb{S}}  R_g(z)(X_g(z)+c)g(z) dz  -\mu \epsilon^{\gamma^2/2}  \int_{\mathbb{S}} e^{\gamma \phi_{\epsilon}(z)}  dz  }    ]   dc   \\
 & =     \int_\R   \E  [  e^{(\sum_i \alpha_i-2Q) c }    \prod_{i=1}^n \epsilon^{\alpha_i^2/2} e^{\alpha_i (X_{g,\epsilon}(z_i) + \frac{Q}{2} \ln g(z_i))}  e^{  -\mu \epsilon^{\gamma^2/2}e^{\gamma c}  \int_{\mathbb{S}} e^{\gamma (X_{g,\epsilon}(z) + \frac{Q}{2} \ln g(z))}  dz  }    ]   dc   \\
  \end{align*} 
   Now, the first step is to get rid of the vertex  fields $ V_{\alpha_i,\epsilon}(z_i) $ in the above expression since they do not converge pointwise as $\epsilon$ goes to $0$. First, we have by \eqref{devvar} that 
   \begin{equation}\label{variancecomput}
   \E[  (  \sum_{i=1}^n \alpha_i X_{g,\epsilon}(z_i)  )^2 ]=  ( \sum_{i=1}^n \alpha_i^2 ) \ln \frac{1}{\epsilon}-\frac{1}{2} \sum_{i=1}^n \alpha_i^2 \ln g(z_i)  + \sum_{i \not = j} \alpha_i \alpha_j G_g(z_i,z_j) + ( \sum_{i=1}^n \alpha_i^2 ) C+o(1)
   \end{equation}
   where $o(1)$ converges to $0$ when $\epsilon$ goes to $0$. We set
    \begin{equation*}
   Y_\epsilon= \sum_{i=1}^n \alpha_i X_{g,\epsilon}(z_i). 
   \end{equation*}
 If we apply the Girsanov theorem with the variable $Y_\epsilon$ and the field $X_{g,\epsilon}(z)$, we get using \eqref{variancecomput} that up to $e^{O(1)}$ terms we have
    \begin{align*}
   &  < \prod_{i=1}^n V_{\alpha_i,\epsilon}(z_i) > /    Z_{\mathrm{GFF}}(g) \\
     & =    e^{\frac{1}{2}\sum_{i \not = j} \alpha_i \alpha_j G_g(z_i,z_j)  }  \prod_{i=1}^n g(z_i)^{\frac{\alpha_i Q}{2}-\frac{\alpha_i^2}{4}} \int_\R  \E  [  e^{(\sum_i \alpha_i-2Q) c }     e^{  -\mu \epsilon^{\gamma^2/2}e^{\gamma c}  \int_{\mathbb{S}} e^{\gamma (X_{g,\epsilon}(z) + H_{(z_i,\alpha_i),\epsilon}(z) + \frac{Q}{2} \ln g(z))}  dz  }    ]   dc,   \\
  \end{align*} 
where $H_{(z_i,\alpha_i),\epsilon}(z)= \gamma \sum_{i=1}^n \alpha_i G_{g,\epsilon}(z_i,z)$ with $G_{g,\epsilon}(z,y)= \E[ X_{g,\epsilon}(z)X_{g,\epsilon}(y)  ] $. We set
   \begin{equation*}
  Z_{(z_i,\alpha_i),\epsilon}(dz)= e^{\gamma \sum_{i=1}^n \alpha_i G_{g,\epsilon}(z_i,z)} M_{\gamma,\epsilon} (dz),
  \end{equation*} 
where $M_{\gamma,\epsilon} (dz)= e^{\gamma (X_{g,\epsilon}(z)  + \frac{Q}{2} \ln g(z))}dz$. Now, we make the change of variables $$u=\epsilon^{\gamma^2/2}e^{\gamma c} Z_{(z_i,\alpha_i),\epsilon}(\S)  $$ in the above formula which leads to 
     \begin{align*}
   &  < \prod_{i=1}^n V_{\alpha_i,\epsilon}(z_i) > /   Z_{\mathrm{GFF}}(g) \\
     & =  \frac{1}{\gamma}  e^{\frac{1}{2}\sum_{i \not = j} \alpha_i \alpha_j G_g(z_i,z_j)  } \:  \prod_{i=1}^n g(z_i)^{\frac{\alpha_i Q}{2}-\frac{\alpha_i^2}{4}} \: \Gamma (\frac{\sum_i \alpha_i-2Q}{\gamma} ,\mu) \:   \E  [  (  Z_{(z_i,\alpha_i),\epsilon}(\S)^{-\frac{\sum_i \alpha_i-2Q}{\gamma} }   ] .   \\
  \end{align*} 
 In particular, since $G_{g,\epsilon}(z,y)$ converges pointwise to $G_g(z,y)$ for $z \not = y$, it is natural to expect  in view of lemma \ref{firstseib} that $ \E  [  (  Z_{(z_i,\alpha_i),\epsilon}(\S)^{-\frac{\sum_i \alpha_i-2Q}{\gamma} }   ]$ converges to $ \E  [  (  Z_{(z_i,\alpha_i)}(\S)^{-\frac{\sum_i \alpha_i-2Q}{\gamma} }   ] $ as $\epsilon$ goes to $0$ (we do not prove this here): if we admit this convergence, we get \eqref{expcorrels}.   
   \qed

\subsection{Properties of the theory}
  
  Now, we state that the vertex operators are indeed primary local fields (these relations are called the KPZ relations after Knizhnik-Polyakov-Zamolodchikov \cite{KPZ})

  \begin{proposition}  [KPZ relation, \cite{DKRV}]
  If $\psi$ is a M\"obius transform, we have
   \begin{equation*}
   < \prod_{i=1}^n V_{\alpha_i}(\psi(z_i)) >=  \prod_{i=1}^n  |\psi'(z_i)|^{-2 \Delta_{\alpha_i}} \: < \prod_{i=1}^n V_{\alpha_i}(z_i) >
  \end{equation*}
  where $ \Delta_{\alpha_i}= \frac{\alpha_i}{2}(Q- \frac{\alpha_i}{2})$. 
  \end{proposition}
  Hence, in CFT language, the vertex operators $V_{\alpha}$ are primary local fields with conformal weight $\frac{\alpha_i}{2}(Q- \frac{\alpha_i}{2})$. Therefore, in LQFT, there is an infinite number of primary local fields hence it is a very rich theory. The above KPZ relation, which is an exact conformal covariance statement, should not be confused with the geometric KPZ relations proved in Duplantier-Sheffield \cite{cf:DuSh} and Rhodes-Vargas \cite{RV1} for the GMC measures defined in theorem \ref{th:existencechaos}. In particular, these geometric formulations of KPZ are very general and do not rely specifically on conformal invariance: they are valid in all dimensions and for all GMC measures defined in theorem \ref{th:existencechaos}.

  Finally, as is common in CFT, one would like to understand the background metric dependence of the theory and express it in terms of  the central charge. More specifically, if $\varphi$ is a smooth bounded function on $\S$, we can consider the metric $e^{\varphi(z)} g(z)|dz|^2$. Then all the formulas of Riemannian geometry of subsection \ref{sub:geom} are valid in this new metric by replacing the function $g(z)$ by the function $e^{\varphi(z)} g(z)$. One can also define a GFF with vanishing mean $X_{e^{\varphi} g}$ in this new metric, etc... Therefore, one can similarly define correlations $< \prod_{i=1}^n V_{\alpha_i}(z_i) >_{e^{\varphi} g}$ by formula \eqref{defcorrels} where one replaces $g$ with the metric $e^{\varphi} g$. The relation between the two correlation functions is given by the so-called Weyl anomaly formula:  
  
  \begin{proposition} [Weyl anomaly, \cite{DKRV}]\label{Weylanom}
  If $\varphi$ is a smooth bounded function on $\S$, we have 
  \begin{equation}\label{Weylanomaly}
  < \prod_{i=1}^n V_{\alpha_i}(z_i) >_{e^{\varphi} g}=  e^{\frac{c_L}{96\pi} \int_{\S} ( |\nabla_g \varphi|^2(z)+ 2R_{g} (z) \varphi(z)) \, g(z) dz}< \prod_{i=1}^n V_{\alpha_i}(z_i) >
  \end{equation}
  where $c_L= 1+6 Q^2$. Hence LQFT is a CFT with central charge $c_L$.
  \end{proposition}
  
  In CFT, the above property can be seen as a definition of the central charge. There are other ways to see the central charge of the model but we will not present them here. Since the function $\gamma \mapsto 1+6 (\frac{\gamma}{2}+\frac{2}{\gamma})^2$ is a bijection from $]0,2[$ to $]25,\infty[$, the Weyl anomaly formula \eqref{Weylanomaly} shows that LQFT can be seen as  a family of CFTs with central charge varying continuously in the range $]25, \infty[$. Hence, LQFT is an interesting laboratory to check rigorously the general CFT formalism developped in physics following the seminal work of Belavin-Polyakov-Zamolodchikov \cite{BPZ}; LQFT should also arise as the scaling limit of many models in statistical physics (just like the SLE introduced by Schramm \cite{Sch} which is a family of continuous random curves corresponding to a geometrical construction of CFTs with central charge ranging continuously in $]-\infty,1]$).

\subsection{The Liouville measures}
  
 As mentioned in subsection \ref{sub:introCFT}, one can usually (but not always) define primary local fields as random distributions. In the context of LQFT, one can indeed  construct the vertex operators $V_\alpha(z)$ as random distributions in the sense of Schwartz; in fact, since the approximate vertex operators \eqref{approxvertex} are positive random functions, one can in fact show that they converge in the space of random measures hence  $V_\alpha(z)$ can be defined as random measures. Of particular interest is the case $\alpha=\gamma$ on which we will focus in this subsection. To be more precise, let us fix $n$ points $z_i$ with $n \geq 3$. We want to define the random measure $V_\gamma(z) dz$ under the formal probability measure $F\mapsto< F \prod_{i=1}^n V_{\alpha_i}(z_i) >/<\prod_{i=1}^n V_{\alpha_i}(z_i)>$. In this context, we denote the underlying probability space $\E^{(z_i,\alpha_i)}[ . ]$. In view of the definition \eqref{defcorrels}, this leads to the following definition of the Liouville measure (where one just inserts a functional of the measure in the correlation function): if $F$ is a functional defined on measures we have
  \begin{align*}
  & \E^{(z_i,\alpha_i)}_\mu[  F( V_\gamma(z) dz   ) ]  \\
  & = Z_{GFF}(g)  \:   \underset{\epsilon \to 0 }{ \lim}   \int_\R \E  [ F(  V_{\gamma,\epsilon}(z) dz  )  \prod_{i=1}^n   V_{\alpha_i,\epsilon}(z_i)  e^{ - \frac{Q}{4 \pi} \int_{\mathbb{S}}  R_g(z)(X_g(z)+c)  -\mu \epsilon^{\gamma^2/2}  \int_{\mathbb{S}} e^{\gamma \phi_{\epsilon}(z)}  dz  }    ]   dc  /  <\prod_{i=1}^n V_{\alpha_i}(z_i)>  \\
  \end{align*}    
  
 Like for the correlation functions, we can obtain a very explicit expression for these Liouville measures in terms of GMC measures. Along the same line as the proof of the correlations, one can show the following explicit expression for the Liouville measure (with the notations of proposition \ref{seibergprop})
  \begin{equation}\label{Liouville}
   \E^{(z_i,\alpha_i)}_\mu[  F( V_\gamma(z) dz   ) ]  =  \frac{ \E \left [    F(  \xi  \frac{Z_{(z_i,\alpha_i)}(dz)}{ Z_{(z_i,\alpha_i)}(\S) }  ) Z_{(z_i,\alpha_i)}(\S)^{-\frac{\sum_i \alpha_i-2Q}{\gamma}} \right  ]   }{  \E \left [     Z_{(z_i,\alpha_i)}(\S)^{-\frac{\sum_i \alpha_i-2Q}{\gamma}}  \right ]  }
  \end{equation}  
  where $\xi$ is an independent variable with density the standard $\Gamma$-law density $\frac{1}{Z}e^{- \mu x} x^{\frac{\sum_i \alpha_i-2Q}{\gamma}-1} dx$ on $\R_{+}$ (where $Z$ is a normalisation constant to make the integral of mass $1$). We can get rid of the $\xi$ variable by conditioning the measure to have volume $1$. This leads to the unit volume Liouville measures we will denote $V_{\gamma}^1(z) dz$:
    \begin{equation}\label{Liouville_unitvolume}
   \E^{(z_i,\alpha_i)}[  F( V_{\gamma}^1(z) dz   )]  =  \frac{ \E \left [    F(   \frac{Z_{(z_i,\alpha_i)}(dz)}{ Z_{(z_i,\alpha_i)}(\S) }  ) Z_{(z_i,\alpha_i)}(\S)^{-\frac{\sum_i \alpha_i-2Q}{\gamma}} \right  ]   }{  \E \left [     Z_{(z_i,\alpha_i)}(\S)^{-\frac{\sum_i \alpha_i-2Q}{\gamma}}  \right ]  }
  \end{equation}  
  
  One can notice that the $\mu$ dependence has disappeared in the expression of the unit volume Liouville measure. However, the unit volume Liouville measure is not a specific GMC measure (divided by its total mass to have volume $1$) as there is still the $Z_{(z_i,\alpha_i)}(\S)^{-\frac{\sum_i \alpha_i-2Q}{\gamma}}$ term in expression \eqref{Liouville_unitvolume}: this term really comes from the interaction term \eqref{interacterm} in the Liouville action \eqref{Liouvilleaction}. Though the Liouville measures are defined when the $(\alpha_i)_{1 \leq i \leq n}$ satisfy the Seiberg bounds \eqref{Seiberg}, one can show that the unit volume measures exist under the less restrictive conditions     
  \begin{equation}\label{Seiberg_soft}
  \forall i, \; \alpha_i<Q \;  \; \;  \; \text{and}     \; \; \; \; Q-\frac{\sum_{i=1}^n \alpha_i}{2}< \frac{2}{\gamma}  \wedge \min_{1 \leq i \leq n} (Q-\alpha_i),
 \end{equation}
  where $x  \wedge y$ denotes the minimum of $x$ and $y$.

  Among the unit volume Liouville measures, one has a very special importance in relation to planar maps: the one where $n=3$ and for all $i$ we have $\alpha_i=\gamma$ (one can check that for all $\gamma$ in $]0,2[$, this choice of $(\alpha_i)_{1 \leq i \leq n}$ satisfies \eqref{Seiberg_soft}). By conformal invariance, we can consider the case $z_1=0$, $z_2=1$ and $z_3=\infty$. In this case, the measure has a very special conformal invariance conjectured on the limit of planar maps called invariance by rerooting. In words, if you sample a point $x$ according to the measure and send $0$ to $0$, the point $x$ to $1$ and $\infty$ to $\infty$ by a M\"obius transform then the image of the measure by the map has same distribution as the initial measure. More precisely, for a point $x$ different from $0$ and $\infty$ let $\psi_x(z)= z/x$ be the unique M\"obius transform of $\S$ which sends $0$ to $0$, the point $x$ to $1$ and $\infty$ to $\infty$. Then we have the following equality for any functional $F$ defined on measures\footnote{A simple and elegant proof of this property was communicated to us by Julien Dub\'edat.}   
 \begin{equation}\label{invrerooting}
  \E^{ (0,\gamma), (1,\gamma), (\infty,\gamma)   } \left [     \int_\S F( ( V_{\gamma}^1(z) dz  )   \circ \psi_x^{-1}  )   V_{\gamma}^1(x) dx   \right  ] =      \E^{ (0,\gamma), (1,\gamma), (\infty,\gamma)  }[    F(  V_{\gamma}^1(z) dz  )] 
 \end{equation}
  where if $\nu$ is a measure on $\S$ and $f: \S \to \S$ some function, the measure $\nu \circ f^{-1}$ is defined by $(\nu \circ f^{-1} ) (A)= \nu (f^{-1} (A))$ for all Borel sets $A$.

 Finally, we mention that a variant to LQG was developped in a series of works by Duplantier-Miller-Sheffield: see \cite{DMS} and \cite{Sheff}. The framework of these works is a bit different than the one we consider in these notes. Duplantier-Miller-Sheffield consider a GFF version of LQG with no cosmological constant $\mu$ and in particular no correlation functions. In this approach based on a coupling between the GFF and SLE, they construct equivalence classes of random measures (called quantum cones, spheres, etc...) with two marked points and coupled to space filling variants of SLE curves. In some sense, their framework is complementary with the one of \cite{DKRV} which considers random measures with 3 or more marked points. The framework of Duplantier-Miller-Sheffield \cite{DMS} is interesting because it establishes non trivial links between (decorated) random planar maps and the so-called quantum cones, spheres, etc...

\subsection{Conjectured relation with planar maps}

\begin{figure}[h]
\centering
\includegraphics[width=0.7\linewidth]{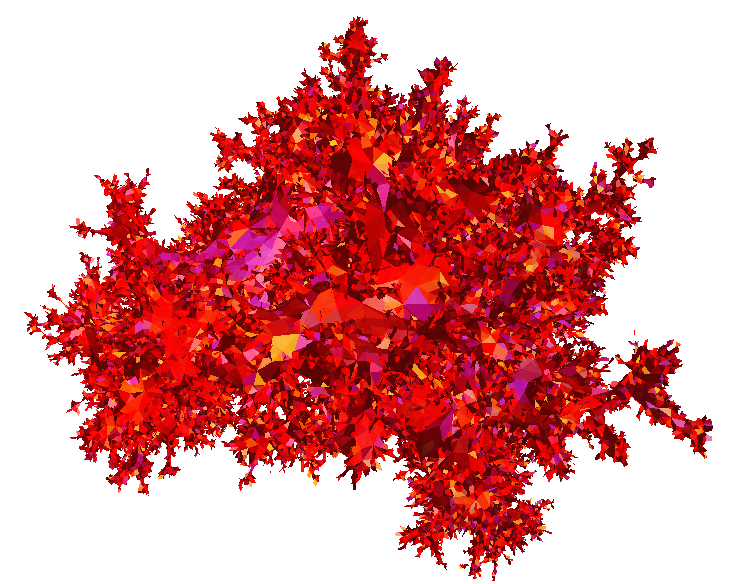}
\caption{Random triangulation with $10^5$ faces (no isometric embedding into the space). Courtesy of F. David}\label{simuDavid}
\end{figure}

Following Polyakov's work \cite{Pol}, it was soon acknowledged by physicists that one should recover LQG as some kind of discretized 2d quantum gravity given by finite triangulations of size $N$ as $N$ goes to infinity (see for example the classical textbook from physics \cite{Amb} for a review on this problem). From now on, we assume that the reader is familiar with the definition of a triangulation of the sphere equipped with a conformal structure: otherwise, he can have a look at the appendix where we gathered the required background. More precisely, let $\mathcal{T}_{N}$ be the set of triangulations of $\S$  with $N$ faces  and  $\mathcal{T}_{N,3}$ be the set of triangulations with $N$ faces and $3$ marked faces (see figure \ref{simuDavid} for a simulation of a random triangulation with $N=10^5$ and sampled according to the uniform measure on $\mathcal{T}_{N}$). We will choose a point in each each marked face: these points are called roots. We equip $T \in \mathcal{T}_{N}$ with a standard conformal structure where each triangle is given volume $1/N$ (see the appendix). The uniformization theorem tells us that we can then conformally map the triangulation onto the sphere $\S$ and the conformal map is unique if we demand the map to send the three roots to prescribed points $z_1, z_2 ,z_3\in\S$. Concretely, the uniformization provides for each face $t\in T$  a conformal map $\psi_t:  t \to \S$ where  $t$ is an equilateral triangle
of volume $\frac{1}{N}$. Then, we denote by $\nu_{T,N}$ the corresponding deterministic measure on $\S$ where $\nu_{T,N}(dz)=|(\psi_t^{-1})'|^2dz$ on each distorted triangle $\tilde{t}$ image of a triangle $t$ by $\psi_t$.  In particular, the volume of the total space $\S $ is $N \times \frac{1}{N}=1$. Now, we consider the random measure $\nu_{N}$ defined by 
\begin{equation}\label{deffinitemap}
\E^{N}[  F( \nu_{N} )  ]= \frac{1}{Z_{N}}  \sum_{T \in \mathcal{T}_{N,3}} F(\nu_{T,N} ),
\end{equation} 
 for positive bounded functions $F$ where $Z_N$ is a normalization constant given by $\# \mathcal{T}_{N,3}$ (the cardinal of the set $\mathcal{T}_{N,3}$). We denote by $\P^{N}$ the probability law associated to $\E^{N}$.

\medskip
We can now state a precise mathematical conjecture:

\begin{conjecture}\label{conjcartes}
Under $\P^{N}$, the family of random measures  $(\nu_{N})_{N \geq 1} $ converges in law as $N \to \infty$ in the space of Radon measures equipped with the topology of weak convergence towards the law of the unit volume Liouville measure given by  \eqref{Liouville_unitvolume} with parameter $\gamma=\sqrt{\frac{8}{3}}$, where $n=3$ and $(z_i,\alpha_i)=(z_i,\gamma)$.
\end{conjecture}

Though such a precise conjecture was first stated in \cite{DKRV}, it is fair to say that such a conjecture is just a clean mathematical formulation of the link between discrete gravity and LQG understood in the 80's by physicists. As of today, conjecture \ref{conjcartes} is still completely open (though partial progress has been made on a closely related question in a paper by Curien \cite{Curien}). One should also mention that a weaker and less explicit variant of conjecture \ref{conjcartes} appears in Sheffield's paper \cite{Sheff}. More precisely, Sheffield proposed a limiting procedure involving the GFF to define a candidate measure for the limit of $(\nu_{N})_{N \geq 1} $ as $N \to \infty$ (see the introduction of section 6 and conjecture 1.(a)); however, he left open the question of convergence of this limiting procedure. Recently, Aru-Huang-Sun \cite{AHS} proved that the limiting procedure does converge and that the limit is the unit volume Liouville measure given by  \eqref{Liouville_unitvolume} with parameter $\gamma=\sqrt{\frac{8}{3}}$, where $n=3$ and $(z_i,\alpha_i)=(z_i,\gamma)$.

Let us consider the case $z_1=0$, $z_2=1$ and $z_3=\infty$ (by conformal invariance, this is no restriction). In this case, one could also consider triangulations with a fourth marked point and send the fourth marked point to $z_3$ in place of the third. Of course, this should not change the limit measure and therefore the limit measure should satisfy the invariance by rerooting property \eqref{invrerooting}.

\begin{figure}[h]
\centering
\subfloat[Circle packing of a triangulation]{\includegraphics[width=0.42\linewidth]{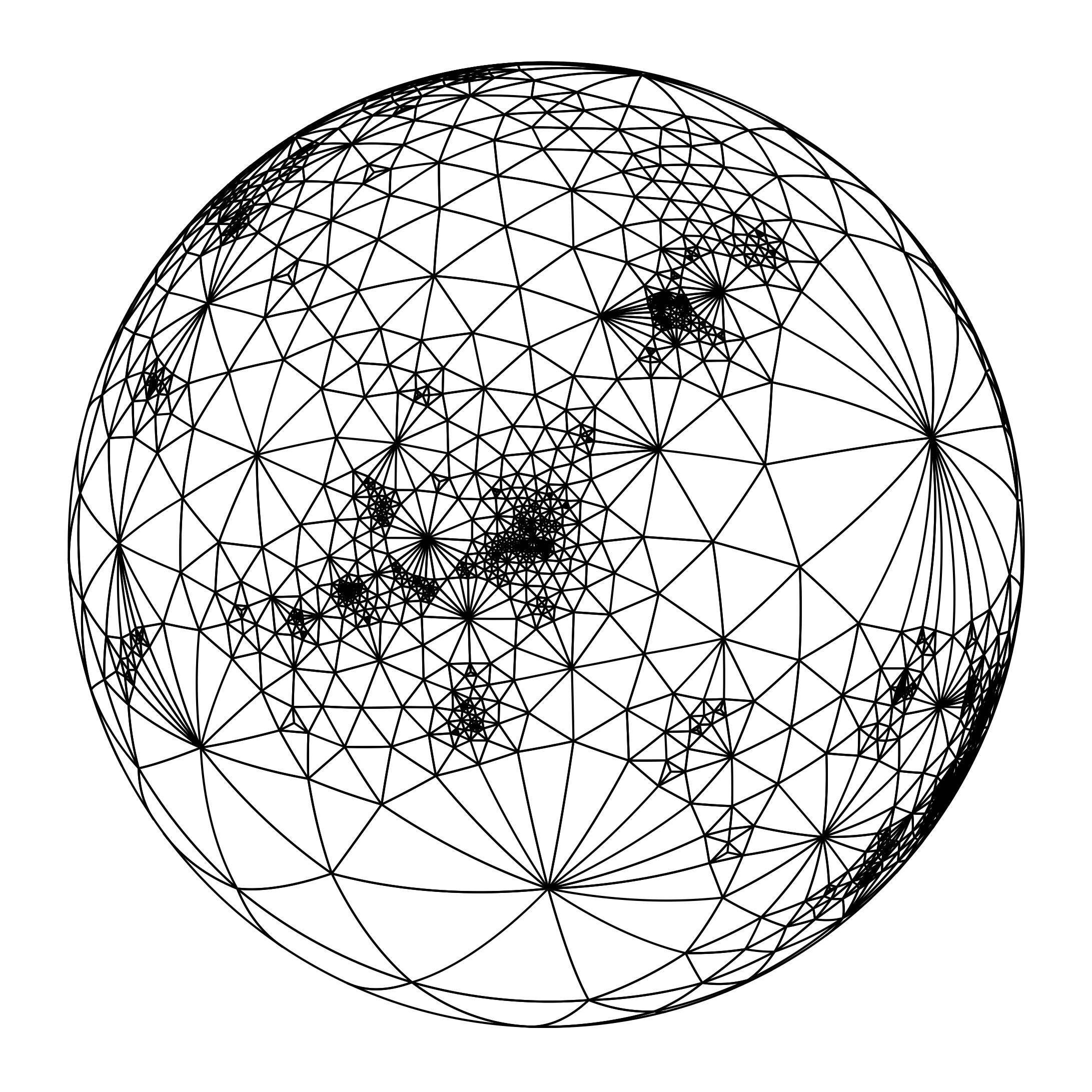}}
\,\,\subfloat[corresponding adjacency circles]{\includegraphics[width=0.42\linewidth]{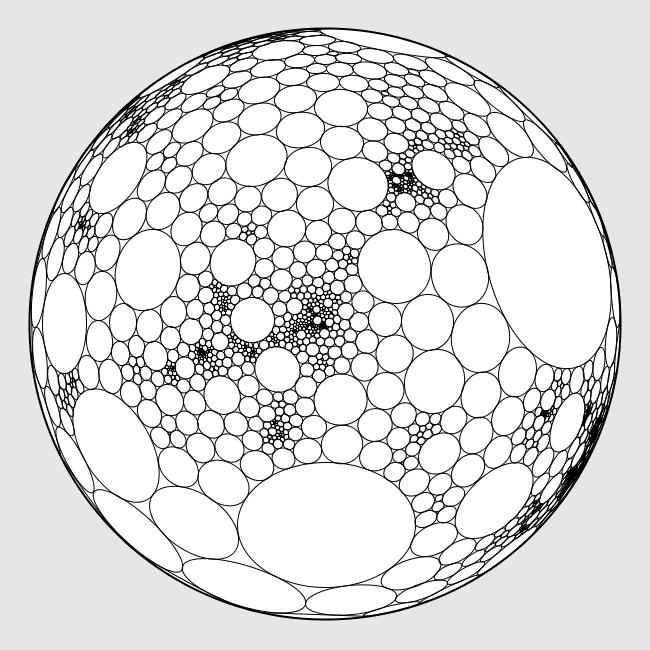}}
\caption{Courtesy of F. David} \label{figcircle}
\end{figure}

Finally, we could also state many variants of conjecture \ref{conjcartes} as it is expected that some form of universality should hold. More precisely, conjecture \ref{conjcartes} should not really depend on the details to define the measure $\nu_{N}$ in \eqref{deffinitemap}. For instance, one expects the same conjecture to hold where $\nu_{T,N}$ could be defined by putting uniform volume $1/N$ in each triangle of the circle packed triangulation: see figure \ref{figcircle} for a circle packed triangulation with large $N$ (however, in this situation, there is a subtelty in the way one fixes the circle packing in a unique way: indeed, M\"obius transforms send circle packings to circle packings but the centers of the circles of the latter are not necessarily the image of the centers of the former by the M\"obius transforms).    

\subsection{On the Ising model at critical temperature}
In this section, we give an account on the recent breakthroughs which occured in the understanding of the Ising model in the plane at critical temperature. This will provide the reader with another example of model where CFT can be made rigorous. Let us start with a few notations. 

On the lattice $\Z^2$ and if $x,y$ are in $\Z^2$ we denote $x \sim y$ the standard adjacency relation.   Let $N$ be a positive integer. We consider the box $\Lambda_N= [|-N,N|]^2$ and its frontier $\partial \Lambda_N= \lbrace x \in {}^c\Lambda_N, \:  \exists y \in \Lambda_N, \: x \sim y \rbrace$. The state space of the model is $\lbrace -1,1 \rbrace^{\Lambda_N}$ and the energy of a spin configuration is given by $$H_N^+(\sigma )= - \sum_{x \in \Lambda_N, \: x \sim y} \sigma_x \sigma_y$$ where we will consider $+$ boundary conditions, i.e. we set the spins in $\partial \Lambda_N$ equal to $1$.

The Ising model on $\Lambda_N$ is then the Gibbs measure $\mu_N$ on the state space $\lbrace -1,1 \rbrace^\Lambda_N$ where the expectation of a functional $F$ is given by 
\begin{equation*}
\mu_{N,\beta}^+(F(\sigma))=  \frac{1}{Z_{N,\beta}}\sum_{ \sigma \in \lbrace -1,1 \rbrace^{\Lambda_N}  } F(\sigma) e^ { -\beta H_N^+(\sigma )} 
\end{equation*} 
where $\beta>0$ is the inverse temperature of the model and $Z_{N,\beta}$ a normalization constant ensuring that $\mu_{N,\beta}^+$ is a probability measure. The model undergoes a phase transition and the critical temperature is explicitly given by $\beta_c= \frac{1}{2} \ln (1+\sqrt{2})$. One can show that the measure $\mu_{N,\beta_c}^+$  converges as $N$ goes to infinity towards a measure $\mu_{\beta_c}$ defined in the full plane, i.e. with state space $\lbrace -1,1 \rbrace^{\Z^2}$ (one can notice that we have removed the superscript $+$ in the full plane measure; indeed one can show that this limit does not depend on the boundary conditions used to define the approximation measures on $\Lambda_N$). 

\begin{figure}[h]
\centering
\includegraphics[width=0.7\linewidth]{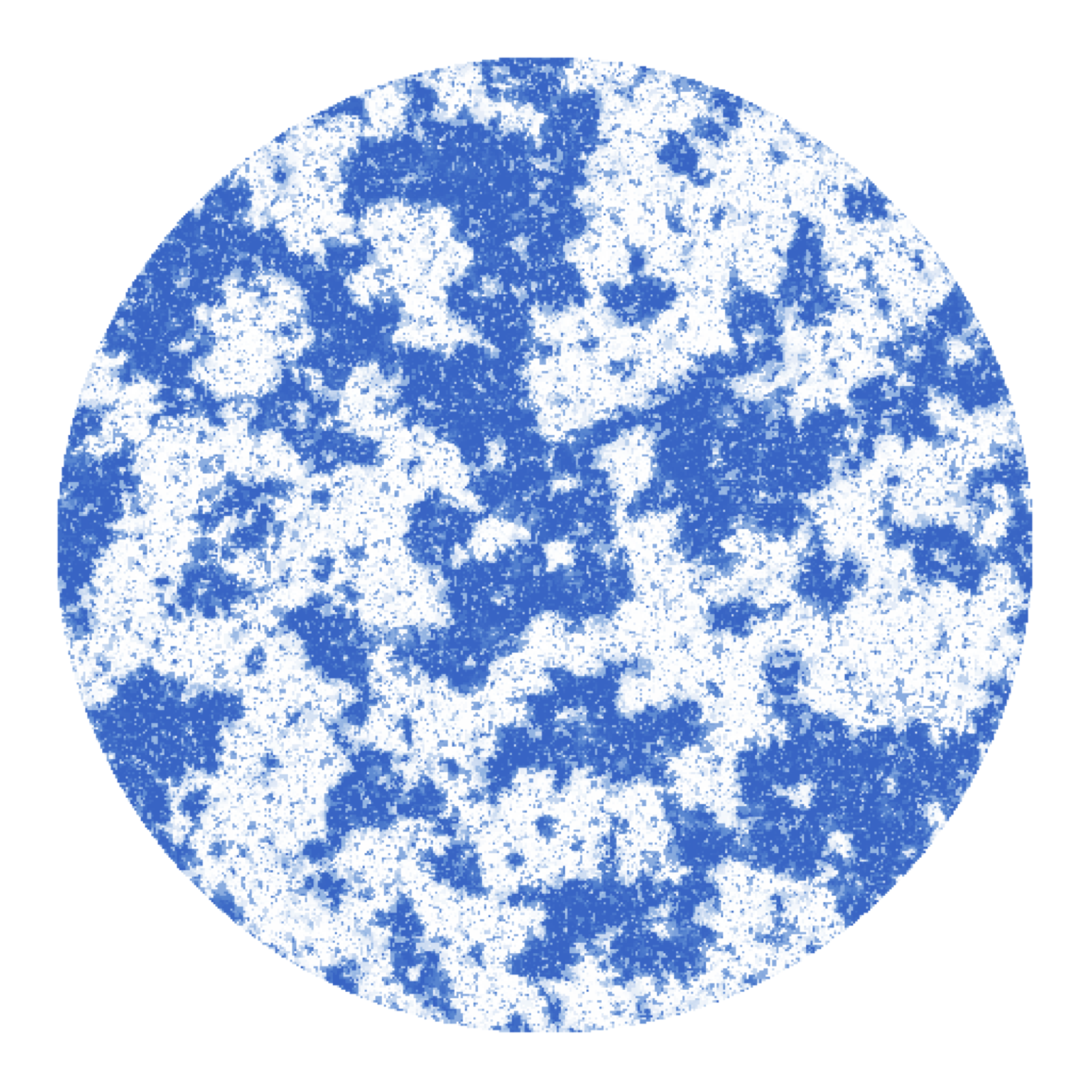}
\caption{Simulated Ising model at critical temperature with free boundary conditions (courtesy of C. Hongler).}\label{figHongler}
\end{figure}

The model was conjectured by physicists to be described by a specific CFT with central charge $c=\frac{1}{2}$ with two primary fields (to be precise there are three primary fields in the theory but the third one is just the constant $1$). We will denote the two primary fields $\sigma(z)$ (the spin field) and $\epsilon(z)$ (the energy density field). We consider the spin field first and set the following definition for non coincident points $z_1, \cdots, z_n$ and $n$ even
\begin{equation*}
<\sigma(z_1) \cdots \sigma(z_n)>:= \left ( 2^{-n/2} \sum_{\mu \in \lbrace -1,1\rbrace^n, \: \sum_i \mu_i=0 }  \prod_{i<j} |z_i-z_j|^{\mu_i \mu_j/2} \right )^{1/2}
\end{equation*}

If $\psi$ is a Mobius transform on the sphere then $|\psi(z)- \psi(y)|= |\psi'(z)|^{1/2} |\psi'(y)|^{1/2} |z-y|$ and therefore
\begin{equation}\label{confinvIsing1}
<\sigma(\psi(z_1)) \cdots \sigma(\psi(z_n))>=  \prod_{i=1}^n |\psi'(z_i)|^{-1/8} <\sigma(z_1) \cdots \sigma(z_n)>
\end{equation}
hence in CFT langage $\sigma$ has conformal weight $\frac{1}{16}$.

Let  $\lfloor .\rfloor$ denote the integer part. For $\epsilon>0$,  we are now interested in the scaling limit as $\epsilon$ goes to $0$ of the discrete spin field $x \mapsto \sigma_{ \lfloor \frac{x}{\epsilon}  \rfloor}$ defined on the rescaled lattice $\epsilon \Z^2$ under the measure $\mu_{\beta_c}$ (see figure \ref{figHongler} for a simulation of the spin field). In view of \eqref{confinvIsing1}, it is natural to rescale the field by the factor $\epsilon^{-1/8}$.

Now the following convergence holds for the rescaled correlations
\begin{equation}\label{convIsing}
\mu_{\beta_c}[  \prod_{i=1}^n  ( \epsilon^{-1/8} \sigma_{ \lfloor \frac{z_i}{\epsilon}  \rfloor} )   ]  \underset{\epsilon \to 0}{\rightarrow} C^n <\sigma(z_1) \cdots \sigma(z_n)>
\end{equation}
where $C$ is a lattice specific constant. This important theorem was proved by Chelkak-Hongler-Izyuorv \cite{CHI} building on the fermionic observable first studied by Smirnov \cite{Smi} and Chelkak-Smirnov \cite{CS}; in fact the main theorem in \cite{CHI} shows the convergence of the rescaled correlations to an explicit expression in any domain (not just the full plane). The convergence result \eqref{convIsing} was also proved independently by Dub\'edat \cite{dubedat} by an exact bosonization procedure (roughly,  bosonization means in this context that there exists an exact relation between the squared correlation functions of the Ising model on a lattice and the correlations of the exponential of the discrete GFF on a lattice).  
As is standard in rigorous CFT, one can define  the limit $\sigma$ as a random distribution. More precisely, Camia-Garban-Newman \cite{CGN} proved that there exists a random distribution $\sigma$ defined on some probability space such that $\epsilon^{-1/8} \sigma_{ \lfloor \frac{x}{\epsilon}  \rfloor}$ converges in law in the space of distributions towards $\sigma$. 

Finally, let us mention that similar results can be proved for the energy density field $\epsilon$. In this case, the properly rescaled (and recentered) energy $\sigma_i \sigma_j$ of a bond between two adjacent vertices $i \sim j$  converges towards the field $\epsilon$ (in the sense of the correlation functions): this is proved in Hongler \cite{Hongler} and Hongler-Smirnov \cite{HS} (in any domain and not just the full plane). It was also proved independently in the full plane by Boutillier and De Tili\`ere on general periodic isoradial graphs \cite{BD,BD1}.  There also exist  explicit formulas for the correlations $<\varepsilon(z_1) \cdots \varepsilon(z_n)>$ of the field $\epsilon$ (but we will not write them here: see \cite{Hongler}) and the field $\epsilon$ has conformal weight $1/2$, i.e. 
\begin{equation}\label{confinvIsing2}
<\varepsilon(\psi(z_1)) \cdots \varepsilon(\psi(z_n))>=  \prod_{i=1}^n |\psi'(z_i)|^{-1} <\varepsilon(z_1) \cdots \varepsilon(z_n)>
\end{equation}
Let us further mention that the energy density field cannot be understood as a random distribution hence $<.>$ is not a real measure in \eqref{confinvIsing2}.

\subsection{Final remarks and conclusion}
In these lecture notes, we introduced the theory of LQFT based on Kahane's GMC theory. More precisely, we introduced the correlation functions and the random measures of the theory. We stated that they satisfy the main assumptions of a CFT on the Riemann sphere. As a comparison and to illustrate the full power of CFT, we also presented in CFT language the recent developments around the Ising model in 2d at the critical point. We would like to stress as a final remark the conceptual difference in the mathematical treatment of the two CFTs. The methods of probabilistic quantum field theory developed in the 1970-1980  around path integral formulations have been up to  now unsuccessful to construct the CFT which describes the scaling limit of the Ising model at critical temperature; it is conjectured that such a construction should exist. Nonetheless, this CFT has been rigorously constructed mathematically by taking the scaling limit of the discrete Ising model hence leaving open the other approach. On the LQFT side, recall that random planar maps (which correspond to discrete gravity) were introduced because defining LQFT by path integral formulations seemed troublesome. The idea was to construct LQFT by taking the scaling limit of large planar maps. However, as we have seen in these lecture notes, a direct construction of LQFT by path integral formulation is feasible whereas proving the convergence of large planar maps is a very difficult topic. Indeed, the convergence has only been established up to now for very specific topologies (of convergence).

\section{Appendix}
  
\subsection{The conformal structure on planar maps}

In this subsection, we recall basic definitions and facts on triangulations equipped with a conformal structure. This part is mostly based on Rhode's paper \cite{GilRho}. A finite triangulation $T$  is a graph you can embed in the sphere such that each inner face has three adjacent edges (the edges do not cross and intersect only at vertices). The triangulation $T$ has size $N$ if it has $N$ faces. We see each triangle $t \in T$ as an equilateral triangle of fixed volume $a^2$ say that we glue topologically according to the edges and the vertices. This defines a topological structure (and even a metric structure). Now, we put a conformal structure on $T$. We need an atlas, i.e. a family of compatible charts. We map the inside of each triangle $t$ to the same triangle in the complex plane. If two triangles are adjacent in the triangulation, we map them to two adjacent equilateral triangles in the complex plane. Now, we need to define an atlas in the neighborhood of a vertex $a$. The vertex $a$ is surrounded by $n$ triangles. We first map these triangles in the complex plane in counterclockwise order and such that each is equilateral. Then we use the map $z \mapsto z^{6/n}$ to "unwind" the triangles (in fact, this unwinds the triangles only if $n>6$) to define a homeomorphism around the vertex $a$. By the uniformization theorem, we can find a conformal map $\psi:T \mapsto \mathbb{C}$ where we send $3$ points in $T$ called roots to fixed points $x_1,x_2,x_3 \in \S$. For each triangle $t$, we can consider $\psi_t$, the restriction of $\psi$ to $t$, as a standard conformal map from $t$ to a distorted triangle $\tilde{t} \subset \S$. It is then natural to equip $\mathbb{C}$ with the standard pullback metric. More precisely, in each triangle $\tilde{t}$ the metric is given by $|(\psi_t^{-1})'(z)|^2 dz$ and then one can define the metric in $\mathbb{C}$ by gluing the metric of each distorted triangle $\tilde{t}$. This metric has conical singularities at the points $\alpha$ of the form $\alpha=\psi(a)$ where $a$ is a vertex of $T$.   

Since $\psi^{-1}$ is analytic, we have $|\psi^{-1}(z)| \approx |z-\alpha|^{n / 6}$ around $\alpha$ (to see this compose $\psi^{-1}$ with the chart $z \mapsto z^{6/n}$ ). Recall that the metric on $\S$ around $\alpha$ is of the form $|(\psi^{-1})'(z)|^2 dz=e^{\lambda(z)} dz$. We have $|(\psi^{-1})'(z)|^2 \approx  |z-\alpha|^{2(n/6-1)}$. Therefore, there is little mass around points $n>6$ and big mass around points $n<6$. This metric has a cone interpretation. If $\theta >0$ is some angle and $C_\theta$ is the corresponding cone, one can put a conformal structure on the cone by the function $\psi: z \mapsto z^{\frac{2\pi}{\theta}}$ in which case the metric is 
\begin{equation*}
|(\psi^{-1})'(z)|^2 dz= \frac{\theta}{2 \pi} |z|^{  2 (\frac{\theta}{2 \pi}-1)  } dz=\frac{\theta}{2 \pi} |z|^{  2 \beta  }  dz
\end{equation*} 
where $\beta= \frac{\theta}{2 \pi}-1 $ is in $]-1, \infty[$. Therefore, around $0$, the average Ricci curvature is then given by
\begin{equation*}
- 2\beta  \int_{|z| \leq 1}  \Delta_z \ln |z| dz= -4 \pi \beta   = 2(2 \pi-\theta)
\end{equation*}  
In the case of triangulations, the angle $\theta$ is related to $n$ by the formula $\theta= \frac{n \pi}{3}$: this means that there is negative curvature (and little mass) around $\alpha$ if $n>6$ and the opposite if $n<6$. 


\hspace{10 cm}

\end{document}